\documentclass{elsart}
\journal{ \ \ }

\renewcommand{\theequation}{\thesection.\arabic{equation}}
\usepackage[centertags]{amsmath}
\usepackage{newlfont}
\usepackage{color}
\usepackage{graphicx}
\usepackage{psfrag}
\usepackage{listings}
\usepackage{color}
\usepackage{booktabs}
\usepackage[T1]{fontenc}
\usepackage[english]{babel}
\usepackage{array}
\usepackage{amsfonts}
\usepackage[utf8]{inputenc}

\definecolor{codegreen}{rgb}{0,0.6,0}
\definecolor{codegray}{rgb}{0.5,0.5,0.5}
\definecolor{codepurple}{rgb}{0.58,0,0.82}
\definecolor{backcolour}{rgb}{0.95,0.95,0.92}

\lstdefinestyle{mystyle}{
    backgroundcolor=\color{backcolour},   
    commentstyle=\color{codegreen},
    keywordstyle=\color{magenta},
    numberstyle=\tiny\color{codegray},
    stringstyle=\color{codepurple},
    basicstyle=\footnotesize,
    breakatwhitespace=false,         
    breaklines=true,                 
    captionpos=b,                    
    keepspaces=true,                 
    numbers=left,                    
    numbersep=5pt,                  
    showspaces=false,                
    showstringspaces=false,
    showtabs=false,                  
    tabsize=2
}
 
\begin{document}
\begin{frontmatter}
\title{Efficient Parameter Selection for Scaled Trust-Region Newton Algorithm in Solving Bound-constrained Nonlinear Systems}

\author{Hengameh Mirhajianmoghadam$^1$,}
\author{ S. Mahmood Ghasemi$^2$}
\address{1. Department of Electrical Engineering, Ferdowsi University of Mashhad, Mashhad, Iran hengame.mirhajian@gmail.com \newline 2.Department of Mathematics, University of Houston, Houston, TX, USA Mahmood@math.uh.edu}
\date{\today}
\begin{abstract}
We investigate the problem of parameter selection for the scaled trust-region Newton (STRN) algorithm in solving bound-constrained nonlinear equations. Numerical experiments were performed on a large number of test problems to find the best value range of parameters that give the least algorithm iterations and function evaluations. Our experiments demonstrate that, in general, there is no best parameter to be chosen and each specific value shows an efficient performance on some problems and a weak performance on other ones. In this research, we report the performance of STRN for various choices of parameters and then suggest the most effective one.

\end{abstract}
\begin{keyword}
Trust-Region Methods, Nonlinear System of Equations, Bound-Constrained Optimizations, Parameter Selection
\end{keyword}
\end{frontmatter}
\section{Introduction}

Optimization plays a key role in contemporary science, including engineering, statistics, computer science, physics, and applied mathematics. If the physics of a phenomenon is properly and comprehensively captured in the associated mathematical model, a huge number of real-world problems translate into solving   mathematical problems; especially in the form of optimization \cite{mo19,layegh18,sch18,layegh20}. Most practical problems may include nonlinear constraints, but constrained optimization remarkably relies on the techniques used in unconstrained optimization. Consequently, the more fruitful unconstrained optimization algorithms are, the more efficient methods can be proposed for the constrained optimization.  Therefore, it is of extreme importance to investigate the properties of the main algorithms used in unconstrained optimization. Here, we examine the performance of a vital algorithm in the trust-region paradigm. We study the efficient parameter selection for the scaled trust-region Newton (STRN) algorithm in solving bound-constrained nonlinear systems. We demonstrate that its performance notably relies on the choices of its hyper-parameters and further suggest how to choose the most effective parameters. 

We consider the STRN algorithm proposed by Bellavia et al.  \cite{bel} and its Matlab solver STRSCNE \cite{bel2}. This numerical algorithm solves the bound-constrained nonlinear system of equations using an affine scaling trust-region method. The problem is
\begin{equation}
F(x)=0; \ \qquad \ x\in \Omega,
\end{equation}
where $F(x)=(F_1(x),\dots, F_n(x))^T$ and $\Omega=\{x\in \mathbb R^n \, | \,  l\leq x \leq u\}.$
The vectors $l\in (\mathbb R\cup - \infty)^n$ and $u\in (\mathbb R\cup  \infty)^n$ are lower and upper bounds, respectively. The function $F$ is continuously differentiable in an open set $X\subset \mathbb  R^n$ containing the n-dimensional box $\Omega$. These kind of systems appear in chemical process modeling and in steady-state simulation \cite{shach2,bull}.\\\\
An approach to solve (1.1) is a bound-constrained nonlinear least square problem:
\begin{align}
\min_{x\in \Omega}f(x):= \frac{1}{2} \left \| F(x) \right \|_2^2
\end{align}
Nonlinear least square problems have been studied in the literature \cite{noc,kanzow,koz,Ulbrich}. Bellavia et al. \cite{bel} generalized the trust-region strategy for unconstrained systems of nonlinear equations to bound-constrained systems and proposed the STRN; a reliable method for tackling (1.1).
This method generates feasible iterates with locally and globally fast convergence properties.
A large set of problems was used to test the efficiency of the STRN. In \cite{bel}, a comparison with the ASTN   \cite{kanzow} and  IGNT  \cite{koz}  and  in \cite{hek,hek3} a comparison between the STRN and NMAdapt \cite{hon} has been performed and superiority of the method has been proved. The method has been widely applied in engineering fields  \cite{fer,man,hos,naj2,pir}.\\
In this paper, we run the iterative algorithm STRN using its implementation in Matlab solver called STRCNE (Scaled Trust-Region Solver for Constrained Nonlinear Equations) on various problems and extract the most useful parameter values. The motivation behind this fine tuning is that the fast-increasing availability of massive data sets has boosted up the use of sophisticated
optimization algorithms \cite{az,naj} with fast convergence. 

The remaining of the paper is as follows: in section 2, we explain the STRN algorithm. Section 3 presents the parameter selecting experiments and the achieved numerical results. Section 4 is the conclusion together with the related figures and tables.

\section{The STRN Algorithm}
Let $x_k\in int(\Omega)$ be the current iteration. Next iteration is  $x_{k+1}=x_k+p_k$ where $p_k$ is computed by solving the following elliptical trust-region subproblem
$$\ \ \ \ \ \ \ \ \ \ \ \ \ \ \ \ \ \ \ \ \ \ \ \ \ \ \ \ \  \ \ \  \min_p\,{m_k(p)\quad \text{ subject to }\quad\left \| D_kp \right \|\leq \Delta_k.} \ \ \ \ \ \ \ \ \ \ \ \ \ \ \ \ \ \ \ \ \ \ \ \ \ \ \  \  (2.1)$$
Here, $\Delta_k$ is  the trust-region size, $D_k=D(x_k)$ is the diagonal scaling matrix such that:
$$D(x)=\text{diag}(|v_1(x)|^{-\frac{1}{2}}, |v_2(x)|^{-\frac{1}{2}},...,|v_n(x)|^{-\frac{1}{2}}),$$
and $v(x)$ denotes the vector function given by:

\[ \left\{ \begin{array}{ll}
         v_i(x)=x_i-u_i & \mbox{ \quad if \quad $( \bigtriangledown f(x) )_i <0,\quad  \text{ and }  \quad  u_i< \infty$};\\
        v_i(x)=x_i-l_i & \mbox{ \quad if \quad $( \bigtriangledown f(x) )_i \ge0,\quad  \text{ and }  \quad  u_i< -\infty$};\\
       v_i(x)=-1  & \mbox{ \quad if \quad $( \bigtriangledown f(x) )_i <0,\quad  \text{ and }  \quad  u_i=\infty$};\\
        v_i(x)=1 & \mbox{ \quad if \quad $( \bigtriangledown f(x) )_i \ge 0,\quad  \text{ and }  \quad  u_i=-\infty$},
        \end{array} \right. \] 
and
$$v=(v_1,...,v_n).$$
The quadratic model for $f$ in (2.1) is as follows:
$$m_k(p)=\frac{1}{2} \left \| F_k'p+F_k \right \|^2=\frac{1}{2}\left \| F_k \right \|^2+F_k^TF_k'p+\frac{1}{2} p^TF_k'^TF_k'p=f_k+\bigtriangledown f_k^Tp+\frac{1}{2}p^TF_k'^TF_k'p.$$
The scaled steepest descent direction $d_k$ is given by
$d_k=-D_k^{-2}\bigtriangledown f_k.$
The trial step $p_k$ can  be calculated by the following procedure (see \cite{bel,bel2} for detail).\\

\textbf{Procedure to calculate the trial step:}\\\\
Let $\bigtriangledown f_k, D_k, $ and $\Delta_k$ be given.\\\\
1. Calculate the Newton step $p_k^N$ by solving $F_k'p_k^N=-F_k$\\\\
2. If $\left \| D_kp_k^N\right \| \leq\Delta_k$ then set $p_k=p_k^N$ and stop.\\
3. Compute $\widetilde{p}_k^u=-\frac{\left \|D_k^{-1} \bigtriangledown f_k \right \|^2}{\left \|F_k'D_k^{-2}\bigtriangledown f_k \right \|^2}D_k^{-1}\bigtriangledown f_k$\\
4. If $\left \| \widetilde{p}_k^u \right \|\geq \Delta_k,$ then set $\widetilde{p}_k=\Delta_kD_k^{-1}\bigtriangledown f_k/\left \|D_k^{-1}\bigtriangledown f_k \right \|$\\
else
set $  \widetilde{p}_k^N=Dp_k^N$ and compute  $\mu$ solving $\left \| \widetilde{p}_k^u +(1-\mu) (\widetilde{p}_k^N-\widetilde{p}_k^u) \right \|^2=\Delta_k^2 $ and set $\widetilde{p}_k=\widetilde{p}_k^u +(1-\mu) (\widetilde{p}_k^N-\widetilde{p}_k^u)$\\\\
5. Let $p_k=D_k^{-1}\widetilde{p}_k$.\\\\
Now to ensure that the next iterate stays within $\Omega$, lets calculate the step size $\lambda(p_k)$ along $p_k$ to the boundary 
$$\lambda(p_k)=\left\{\begin{matrix}
\infty &  \quad \text{ if } \quad  \ \ \Omega = \mathbb R^n \\ 
 \min_{1\leq i\leq n} \Lambda _i& \quad \text{ if } \quad  \ \ \Omega \subset \mathbb R^n,  
\end{matrix}\right.$$
where
$$\Lambda _i=\left\{\begin{matrix}
 \max\big\{\frac{l_i-x_{k_i}}{p_{k_i}},\frac{u_i-x_{k_i}}{p_{k_i}}\big \}& \quad \text{ if } \quad \ p_{k_i}\neq 0 \\
\infty &  \quad \text{ if } \quad  \ p_{k_i}= 0.
\end{matrix}\right.$$
if $\lambda(p_k)>1$ then $x_k+p_k$ is within $\Omega$, if $\lambda(p_k)\leq 1$ then a step back along $p_k$ is necessary. Parameter $\theta$ controls the amount of truncation:
$$x_{k+1}=x_k+\alpha(p_k)$$
$$\alpha(p_k)=\left\{\begin{matrix}
p_k &  \quad \text{ if } \quad  \lambda(p_k)>1\\ 
 \max\{\theta,1-\left \|p_k\right \|  \} \lambda(p_k)p_k&  \quad \text{ if } \quad \lambda(p_k)\le 1.
\end{matrix}\right.$$

$\theta \in (0,1) $ is a fixed constant and it is one of the  parameters that we will find an optimal range for. \\
In order to warranty sufficient reduction, we have to consider the Cauchy point $p_k^c$, the minimizer of $m_k$ along the scaled steepest descent direction $d_k=-D_k^{-2}\bigtriangledown f_k $, at the other hand new iterate should lie within the trust-region so the steepest descent has to satisfy the trust-region bound \cite{noc} 
$$p_k^c=\tau_k d_k=-\tau_k D_k ^{-2} \bigtriangledown f_k,$$ where 
 $$\tau_k=\text{argmin}_{\tau>0}\{ m_k(\tau d_k) : \left \|D_kd_k \right \| \leq \Delta_k \}=\min\Bigg\{ \frac{\left \|D_k^{-1} \bigtriangledown f_k \right \|^2}{\left \|F_k'D_k^{-2}\bigtriangledown f_k \right \|^2},\frac{\Delta_k}{\left \|D_k^{-1} \bigtriangledown f_k \right \|} \Bigg \}.$$
Then we test if the step $\alpha(p_k)$ satisfies the following condition:
$$\rho_k^c(p_k)=\frac{m_k(0)-m_k(\alpha(p_k))}{m_k(0)-m_k(\alpha(p_k^c))}\geq \beta_1$$
where $\beta_1\in(0,1)$ is constant. If this condition holds, we discard $p_k$ and set $p_k=p_k^c$. \\
The agreement between the model $m_k$ and the merit function $f$ can be achieved by testing the following condition:
$$\rho_k^f(p_k)=\frac{f(x_k)-f(x_k+\alpha(p_k))}{m_k(0)-m_k(\alpha(p_k))}\geq \beta_2$$
where $\beta_2 \in(0,1]$. If this condition holds then $x_k+\alpha(p_k)$ is the next iterate. If not then we have to decrease the trust-region size by
$$\Delta_k=\min\{ \alpha_1 \Delta_k,\alpha_2 \left \|D_k \alpha (p_k) \right \|\}$$
for $0<\alpha_1\leq \alpha_2 <1$ and recalculate new step. \\
Finally, in order to accelerate the convergence rate, we should take big steps by increasing the trust-region radius wisely, If the agreement between the function and the model is strong enough then we shouldn't miss the opportunity to take a better improvement. The condition below satisfies this:
$$\rho_k^f(p_k)=\frac{f(x_k)-f(x_k+\alpha(p_k))}{m_k(0)-m_k(\alpha(p_k))}\geq \beta_3,$$
$\beta_3\in (0,1]$ is a constant such that $\beta_2 < \beta_3< 1.$\\
If this condition holds then $$\Delta_{k+1}=\max\{ \Delta_k,\gamma \left\| D_k \alpha(p_k) \right\| \},\ \ \ \ \ \gamma>1$$ otherwise the trust-region radius remains same. 
\\
The STRSCNE algorithm is as follows:\\\\
Algorithm, The Scaled Trust-Region Solver \cite{bel2}\\
- Initialization:\\
Given $x_0 \in int(\Omega), \Delta_0>0, \theta \in (0,1), 0< \alpha_1\leq \alpha_2<1,\beta_1 \in (0, 1], 0 < \beta_2 < \beta_3 < 1.$\\
- For k = 0, 1, . . . do:\\
1. Compute $ F_ k$ .\\
2. Check for convergence.\\
3. Compute the matrix $D _k$ by using (1.5).\\
4. Compute the matrix $F _k'$.\\
5. Compute $p_k ^N$ by solving the linear system $F_k'p_k^N=-F_k$.\\\\
6. Repeat\\
6.1. Compute an approximate solution $p_k$ of (2.1) by using Procedure to calculate the trial step.\\
6.2. Compute $\tau_ k$ by  and the Cauchy point $p _k^ c$.\\
6.3. Compute $\alpha ( p_ k ) $ and $\alpha ( p_ k^ c )$.\\
6.4. Compute $\rho_ k^ c ( p _k )$.\\
6.5. If $\rho_ k^ c ( p_ k ) < \beta_1$ then set $p_ k = p_ k^ c$.\\
6.6. Set  $\Delta_k^*=\Delta_k$ and decrease $\Delta_k$\\
6.7. Compute $\rho_ k ^f( p_ k )$.\\\\
Until $\rho_k^f(p_k)\geq \beta_2$\\
7. Set $x _{k+1} = x_ k + \alpha( p_k), \Delta_k=\Delta_k^*$\\
8. If $\rho_k^f\geq \beta_3$ then\\
set $\Delta_{k+1}=\max\{\delta_k,\gamma \left\| D_k \alpha (p_k)\right\| \}$\\
else set $\Delta_{k+1}=\Delta_k$\\\\\\
The algorithm is convergent when we have:
$\left\| F_{k+1}\right\| \leq 10^{-8}$.\\\\
Failure happens if  \\
1) A maximum number of iterations are performed.\\
2) A maximum number of F-evaluations are performed.\\
3) The trust-region size is reduced below $10^{-8}$ \\
4) The relative change in the function value satisfies
$$\left\| F_{K+1}-F_k \right\|\leq 10^{-14}\left\| F_k \right\|$$
5) The norm of the scaled gradient of the merit function becomes  small:
$$\left\| D_k^{-1} \bigtriangledown f_k \right\| < 10^{-14} $$
6) The scaling matrix $D_  k $ cannot be computed.

\section{Numerical Experiments of Parameter Tuning}
In our experiments, we used Matlab 2018b. For the initialization, we set $\Delta_0 =1,\theta = 0.99995, \alpha_ 1 = 0.25, \alpha_ 2 = 0.5, \beta_ 1 = 0.1,\beta_ 2 = 0.25, \beta_ 3 = 0.75,\gamma =2.$ The tested problems come from NLE library \cite{shach} accessible through: www.polymath-software.com/library. 
The number of the iterations and the number of the function evaluations follow the same behavior, lots of the problems are indifferent to changes in the parameters but some of them are influenced by different values of the parameters.  In this section, we run the STRN for different values of parameters  in order to find the best values with least iteration and function evaluation. In each section we fix all the parameters as indicated above and impose the variation only on the specific one.\\
The fact is that there is no "Best" parameter. A parameter works pretty good for some problems and pretty bad for other problems, so we have to choose the one with overall better performance. Tables 1-6 show the iterations for different starting points. 

\subsection{$\alpha_1$- controls the size reduction of the trust-region}
Among the 45 studied problems, 32 problems are insensitive for different values of $\alpha_1$ while 13  problems are sensitive for at least one of the initial points.
After consecutive runs of the algorithm, we narrowed them down to some candidate values for each parameter. Table 1 shows these values.
There is no best parameter for $\alpha_1$. For example, $\alpha_1 =0.4$ is the best option for solving Foureq1 starting from the third $x_0$ while for Threeq6 starting from the fourth $x_0$, or seveneq2b starting from Third $x_0$,  $\alpha_1=0.3$ is much better. Generally speaking, a value between (0.4,0.5) is the best choice (Fig 1). It seems that a harsh and severe  shrinking of the trust-region is not always the best approach, we prefer to cut the trust-region size at most in half.
\begin{figure}[h]
\caption{ }
\centering
\includegraphics[width=4in]{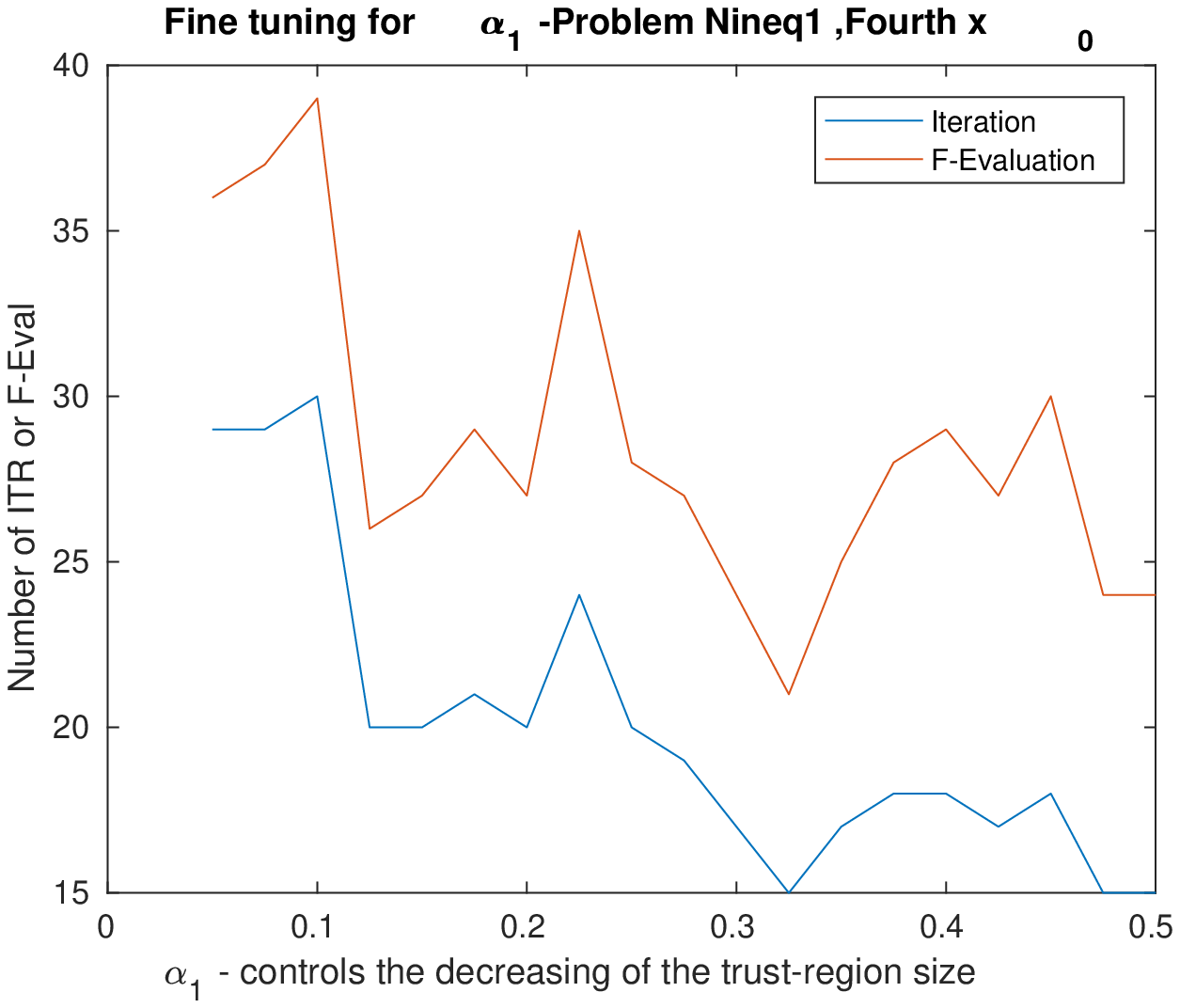}
\end{figure}
\subsection{$\alpha_2$- controls the size reduction of the trust-region }
Among the 45 studied problems,  36 problems are insensitive to different values of $\alpha_2$ while 9  problems are sensitive for at least one of the initial points.
Again, there is no best choice and depending on the problem the performance changes. For example, for Threeq5 (Fig 2) 0.24 and 0.45 are best choices while for other problems they are not suitable.
\begin{figure}[h]
\caption{ }
\centering
\includegraphics[width=4in]{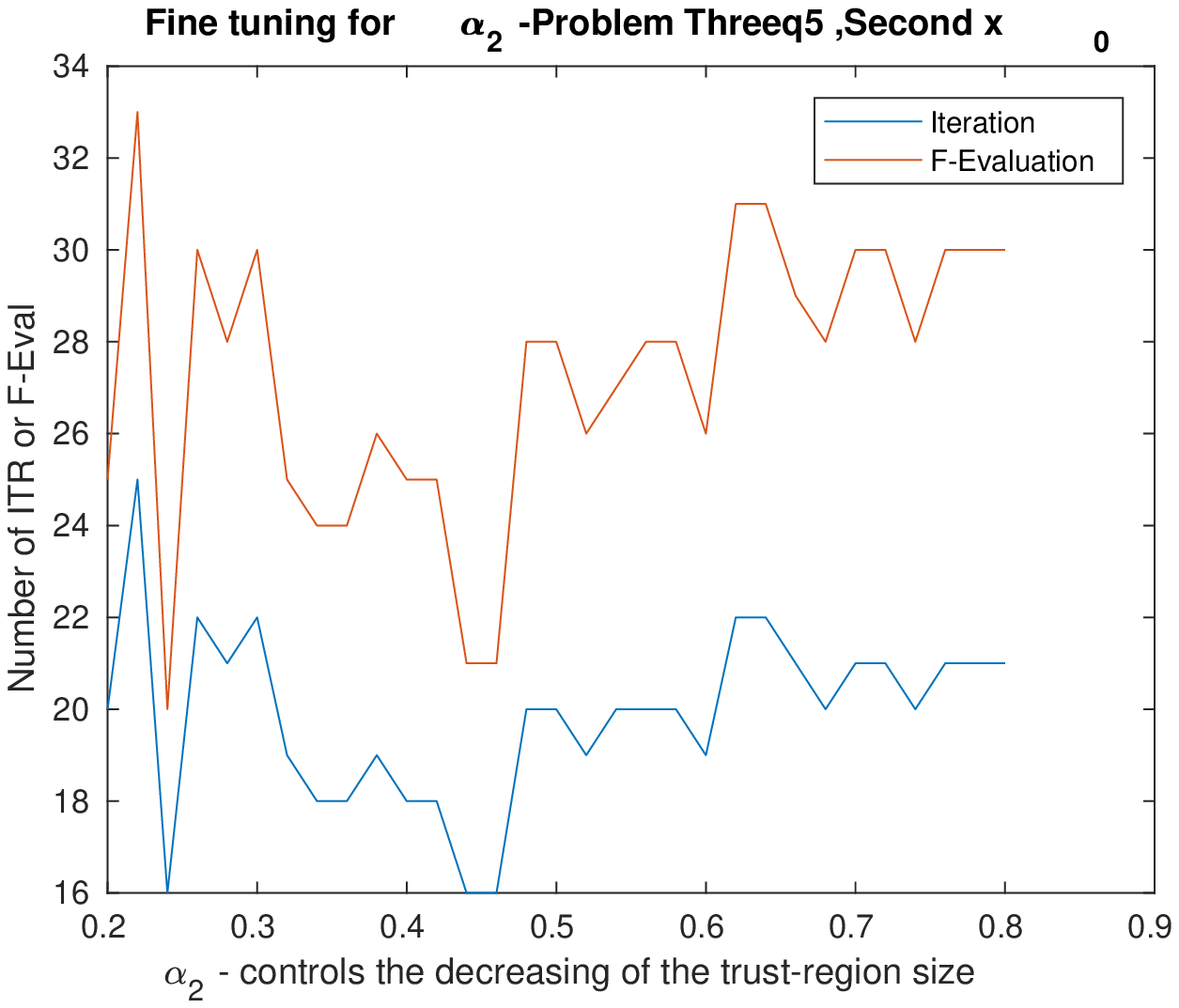}
\end{figure}
by looking at Table 2 the overall performance of $\alpha_2=0.45, 0.6$ are slightly better. $\alpha_2=0.45$ is the best choice for problem Seveneq2b, it converges to the answer in 130 iterations while other selections of the parameter need a large iteration number until getting convergent. $\alpha_2=0.4 $ and $\alpha_2=0.7$ also show pretty good performance.
\subsection{$\beta_1$-used for accuracy requirements }
Among the 45 studied problems,  43 problems are insensitive for different values of $\beta_1$ while 2  problems are sensitive to at least one of the initial points.
By looking at Table 3, as $\beta_1$ increases (The reduction in the model when we use $p_k$ is close to the reduction in the model when we use $p_k^c$) the algorithm becomes faster. $\beta_1=0.2$ shows good performance. The problems are mostly insensitive to the different values of $\beta_1$.
\subsection{$\beta_2$- used to ensure agreement between the model and the objective}
Among  the 45 studied problems,  36 problems are insensitive for different values of $\beta_2$ while 9  problems are sensitive for at least one of the initial points.
From Table 4 clearly $\beta_2=0.15$ is the best choice. As $\beta_2$ increases the agreement between the model and the function becomes harder and the method becomes slower.  We have this situation in problem Nineeq1 (Fig 3):
\begin{figure}[h]
\caption{ }
\centering
\includegraphics[width=4in]{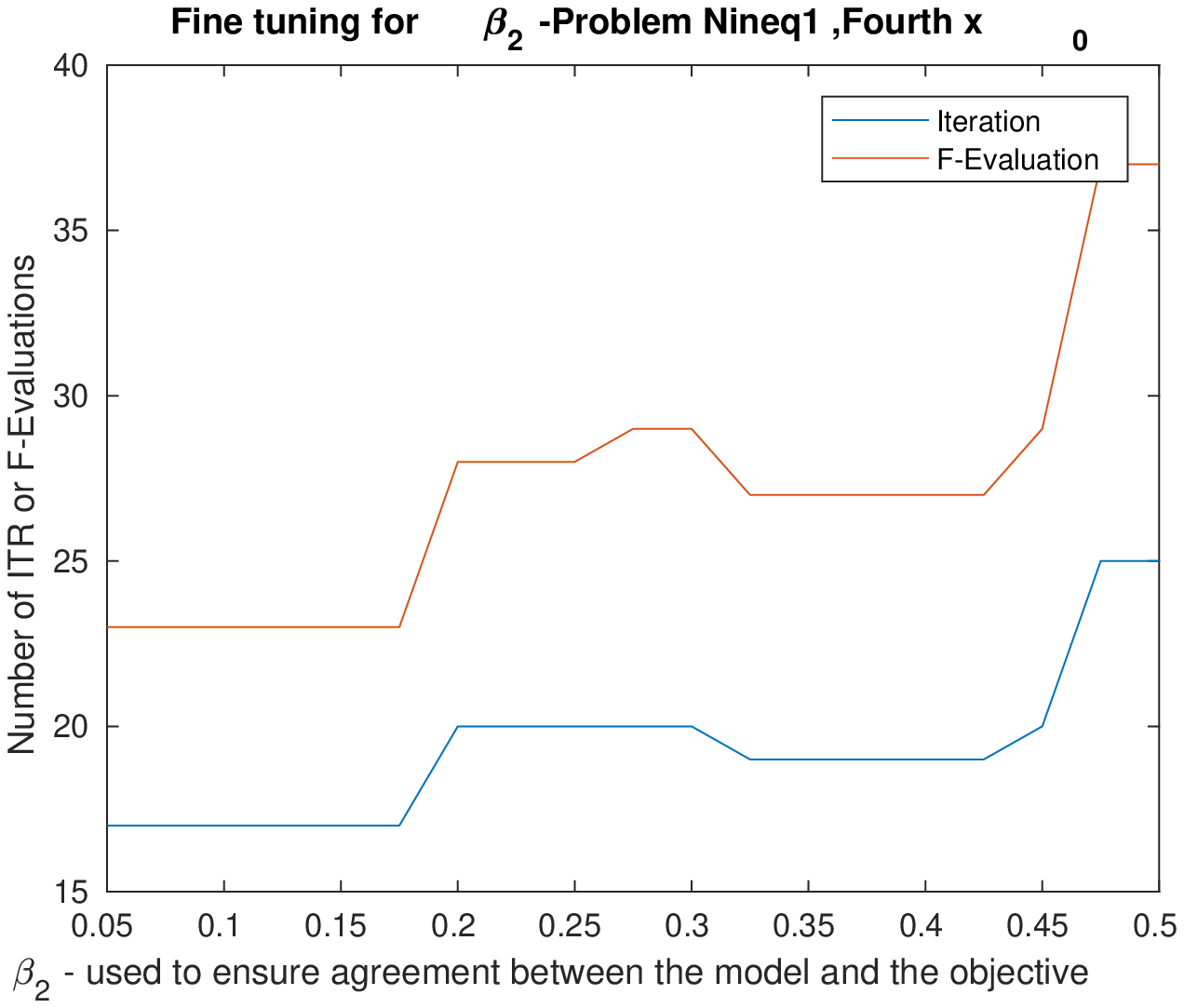}
\end{figure}
\subsection{$\beta_3$-controls the updating of the trust-region size}
Among  the 45 studied problems,  31 problems are insensitive to different values of $\beta_3$ while 14  problems are sensitive for at least one of the initial points.
$\beta_3=0.75,0.6$ are the best choices.  It means that we shouldn't wait for a great agreement between the model and the function and should take the opportunity to increase the trust-region size.

\subsection{$\theta$-used to ensure strictly feasible  iterates }
Among  the 45 studied problems,  28 problems are insensitive to different values of $\theta$ while 17  problems are sensitive for at least one of the initial points.
Looking at Table 6,  $\theta=0.9$ gives the best results. $\theta=0.95 $ is also a good choice, it works as good as the traditional $\theta=0.99995$, so a big truncation of $p_k$ shows as good performance as small truncation. $\theta=0.7$ is working surprisingly good. It solves the problem Threeq4a in the least possible iterations. 
\subsection{$\gamma$-controls the size enlargement  of the trust-region}
Among the 45 studied problems,  21 problems are insensitive to different values of $\gamma$ while 24  problems are sensitive for at least one of the initial points.
Looking at Fig 4 in most of the cases $\gamma=8,10$ show better performance comparing to $\gamma=2$. It seems that taking big steps and going through large trust-region sizes is risky but it worth to take this risk since in most of the situations it shows faster convergence. 
\section{ Conclusion}
The problems Twoeq5a, Twoeq5b, Twoeq7, Threeq5, and 11eq1 are nearly sensitive to changes in the parameters. The sensitive problems that can be affected by different values of the parameters are: Twoeq 6, Threeq1, Threeq4a, Threeq4b, Threeq6, Threeeq6, Foureq1, Sixeq4b, Seveneq2b, Nineq1 and super sensitive problems that can be affected by all of the 7 parameters are Seveneq2b, and Nineq1. \\
The results show that taking big steps by increasing the trust-region size is risky, but it causes faster convergence and it worth to take the risk. Also, we prefer not to cut the trust-region too much, we prefer to stay in the Cauchy direction as much as possible and take easier criteria to check the agreement between the model and the merit function.\\
The interesting point is that for some of the problems (Threeq6, Sixeq4b, Seveneq2b, and 11eq1) only specific values of $\beta_3$ and $\theta$ work well.  We concentrate on $\theta$ since it plays a more important role.  The fact that the method fails for some values of the $\theta$ motivates us to try different values of this parameter before giving up. This can be done by introducing a sequence of parameters as initial values of $\theta$ and not a constant one. That means in STRSCNE algorithm the initialization step should be modified as follows:\\\\
STRSCNE with a variable $\theta$\\
- Initialization:\\
Given $x_0 \in int(\Omega), \Delta_0>0, \theta_i =\{0.7+0.025i\}_{i=0}^{12}, 0< \alpha_1\leq \alpha_2<1,\beta_1 \in (0, 1], 0 < \beta_2 < \beta_3 < 1$, Let $i=0$ and $\theta=\theta_0$.\\\\
Steps 1-8  \\\\
9. If $\left\| F_{k+1}\right\| > 10^{-8}$ then repeat the algorithm. If $\left\| F_{k+1}\right\| \leq 10^{-8}$:\\
  9.1  If Ierr=0 then\\
STOP and report the solution.\\
else set $i=i+1$, $\theta=\theta_i$ and go to step 1 unless $i=12.$\\\\
It is either convergent to the solution or repeats the algorithm for different values of $\theta$ before reporting a defeat.\\
This can be a useful approach if the goal is solving a problem several times in a restricted amount of time. For example Problem Sixeq4b  starting from third $x_0$ and $\theta=0.6$ is convergent to the solution in 9 iterations while for $\theta=0.99995$ the number of iterations is 350. For the first choice of $\theta$ it takes 1.24 seconds and for the second choice, it takes 41.72 seconds. This happens in several other problems too.  
\newline
\newline

\begingroup\makeatletter\def\f@size{10}\check@mathfonts
\begin {table}[!htb]
\centering
\tiny
\caption {The Algorithm's Iterations for the selected values of $\alpha_1$}\label{tab:title}
\scalebox{1}{
\begin{tabular}{ |>{\small}c|c|c|c| c|c|c| } 
 \hline
First $x_0$\\
\hline
Problem &  $\alpha_1=0.2$ & $\alpha_1=0.3$  & $\alpha_1=0.4$ & $\alpha_1=0.5$ & $\alpha_1=0.6$&$\alpha_1=0.7$\\ 
Twoeq7 &  9 & 9 &7&7&7&7 \\ 
 \hline
Second $x_0$\\
\hline
Threeq1 &  42 & 32&32&32&32&32 \\ 
Threeq4b &  7&7& 7&7&7&6\\ 
foureq1 &  10&10& 10&10&16&15\\ 
 \hline
Third $x_0$\\
\hline
twoeq5b &  7&8 & 7&7&7&7 \\ 
Twoeq6 &  10&9& 12&9&9&9\\ 
Threeq1& 31&35 &30&32 &32&32\\
 Threeq4a &  62&54& 118&70&62&62\\ 
Threeq4b&  10&7& 7&7&7&7\\ 
Foureq1& 13 &13&11 &12&12&12\\ 
Sixeq4b &  326 &301 &299&334&334&334 \\ 
Seveneq2b& 117 & 35 &112&84&84&84 \\ 
Seveneq3a &  8&9& 7&7&7&7\\ 
 \hline
Fourth $x_0$\\
\hline
twoeq5a &  9 &7&9&7&7&7 \\ 
Twoeq6 & 6 &7&6&8&8&8\\ 
Threeq4a& 51 &52&53 &54&54&54\\ 
Threeq6 &  100&76& 225&192&192&192\\ 
Nineq1&  20&17& 18&15&15&15\\ 
 \hline
\end{tabular}
}
\end {table}
\endgroup

\begin {table}[!htb]
\centering
\tiny
\caption {The Algorithm's Iterations for the selected values of $\alpha_2$}\label{tab:title}
\scalebox{1}{
\begin{tabular}{ |c|c|c|c| c|c|c| } 
 \hline
Second $x_0$\\
\hline
Problem &  $\alpha_2=0.3$ & $\alpha_2=0.4$  & $\alpha_2=0.45$ & $\alpha_2=0.5$ & $\alpha_2=0.6$&$\alpha_2=0.7$\\ 
Twoeq6 & 14 & 16&16&16&13&17 \\ 
Threeq1 &  45&32& 32&32&32&32\\ 
Threeq5 &  22&18& 16&20&19&21\\ 
foureq1 &  13&13& 10&16&13&12\\ 
 \hline
Third $x_0$\\
\hline

Twoeq6 &  12&13& 14&13&11&10\\ 
Threeq4a &  52&65& 65&67&67&67\\ 
Threeq4b&  9&7& 7&7&7&7\\ 
Sixeq4b &  321 &493 &495&496&480&497 \\ 
Seveneq2b & 277 & 252 &130&199&165&157 \\ 
 \hline
Fourth $x_0$\\
\hline
Twoeq6 & 6 &12&6&6&6&6\\ 
Threeq4a& 51 &51&52 &59&53&54\\ 
Nineq1 &  22&18& 20&20&19&20\\ 
 \hline
\end{tabular}
}
\end {table}

\begin {table}[!htb]
\centering
\tiny
\caption {The Algorithm's Iterations for the selected values of $\beta_1$}\label{tab:title}
\scalebox{1}{
\begin{tabular}{ |c|c|c|c| c|c|c| } 
 \hline
First $x_0$\\
\hline
Problem &  $\beta_1=0.05$ & $\beta_1=0.1$  & $\beta_1=0.15$ & $\beta_1=0.2$ & $\beta_1=0.25$ \\ 
Nineq1 &  200 & 200 &21&21&21 \\ 
 \hline
Third$x_0$\\
\hline

Seveneq2b&199 &  199 & 199&179&179 \\ 
\hline

\end{tabular}
}
\end {table}
\begin {table}[!htb]
\centering
\tiny
\caption {The Algorithm's Iterations for the selected values of $\beta_3$}\label{tab:title}
\scalebox{1}{
\begin{tabular}{ |c|c|c|c| c|c|c| } 
 \hline
First $x_0$\\
\hline
Problem &  $\beta_3=0.6$ & $\beta_3=0.7$  & $\beta_3=0.75$ & $\beta_3=0.8$ & $\beta_3=0.85$\\ 
Twoeq7 &  9 & 10 &8&8&8 \\ 
Threeq1 &  5&5& 5&6&6\\ 
Threeq6&  8&8& 8&8&25\\ 
 \hline
Second $x_0$\\
\hline
Threeq1 &  38 & 40&32&45&58 \\ 
Threeq5 &  20&20& 20&42&93\\ 
Threeq6 &  116&303& 305&307&-\\ 
Threeq7 &  18 & 13 &13&13&13 \\ 
Seveneq2b &  17&14& 14&14&14\\ 
 \hline
Third $x_0$\\
\hline

twoeq5b &  8 & 8&7&7&7 \\ 
Twoeq6&  13&13& 13&13&10\\ 
Threeq1& 38 &31&24 &24&39\\ 
Threeq2 &  39 & 39&39&44&44 \\ 
Threeq4a &  44&53& 67&64&75\\ 
foureq1& 11 &13&13 &13&20\\ 
Sixeq4b&  380 &361 &350&342&337 \\ 
Seveneq2b & 207 & 202 &199&206&178 \\ 
11eq1 &  19&19& 19&19&20\\ 
 \hline
Fourth $x_0$\\
\hline
Threeq4a&  33 &54&59&60&64 \\ 
Threeq6 & 313 &96&95&93&-\\ 
Nineq1& 19 &19&20 &20&20\\ 
 \hline
\end{tabular}
}
\end {table}
\clearpage
\begin {table}[!htb]
\centering
\tiny
\caption {The Algorithm's Iterations for the selected values of $\theta$}\label{tab:title}
\scalebox{1}{
\begin{tabular}{ |c|c|c|c| c|c|c| } 
 \hline
First $x_0$\\
\hline
Problem &  $\theta=0.6$ & $\theta=0.7$  & $\theta=0.8$ & $\theta=0.9$ & $\theta=0.95$ & $\theta=0.99995$\\ 
Nineq1&  14&12& 14&8&51&200\\ 
Teneq1a &  15&14& 14&14&14&14\\ 
 \hline
Second $x_0$\\
\hline

Twoeq6 &  16 & 17&15&14&14&16 \\ 
Threeq1 &  42&32& 40&32&32&32\\ 
Threeq5 &  15&15& 15&15&20&20\\ 
Threeq6&  97&122& 134&102&319&305\\ 
Threeq7 &  11& 11 &10&12&14&13 \\ 
foureq1 &  13&11& 13&10&13&16\\ 
Sixeq4b &  -&10& -&10&9&7\\ 
Teneq1b &  95& 80&97&90&79&14\\ 
Teneq2b &  13&12& 12&11&11&11\\ 
 \hline
Third $x_0$\\
\hline

Twoeq6 &  18 & 15&15&16&12&13 \\ 
Twoeq7&  8&8&8&7&7&7\\ 
Threeq4a &  70&54& 57&58&67&67\\ 
Threeq4b& 7 &6&7 &6&6&7\\ 
foureq1 &  12&12& 13&13&13&13\\ 
Sixeq4b &  9&9& 170&260&253&350\\ 
Seveneq2b& 40 &153&- &21&41&199\\ 
Teneq2a &  11 & 11 &10&9&9&9 \\ 
11eq1 &  -&-& -&-&22&19\\ 
 \hline
Fourth $x_0$\\
\hline
Threeq4a &  47 &62&48&51&52&59 \\ 
Sixeq4b& 9 &8&8 &7&7&7\\ 
Nineq1& 24 &23&22 &22&19&20\\ 
 \hline
\end{tabular}
}
\end {table}
\begin{figure}[h]
\caption{Different $\gamma$ s performance on sensitive problems}
\centering
\includegraphics[width=4in]{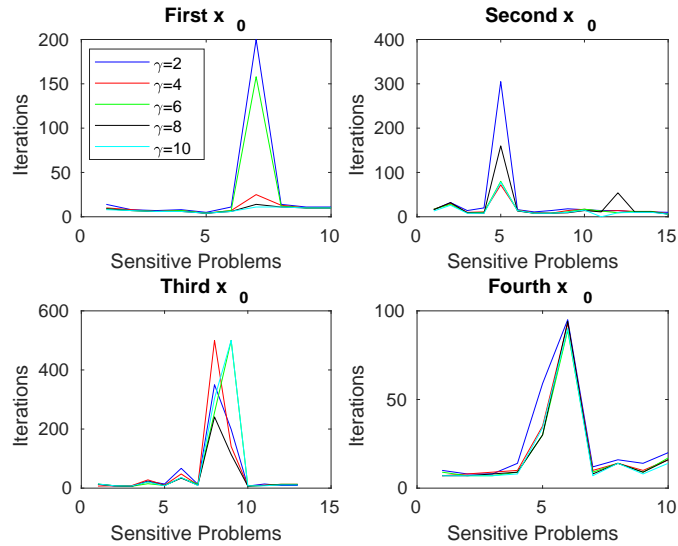}
\end{figure}

\end{document}